\documentclass[12pt]{amsart}

\makeatletter

\theoremstyle{plain}

\newtheorem{thm}{Theorem}[section]

\newtheorem{cor}[thm]{Corollary} 

\newtheorem{lem}[thm]{Lemma} 
\newtheorem{prop}[thm]{Proposition} 
\theoremstyle{definition}
\newtheorem{defn}[thm]{Definition}
\theoremstyle{remark}
\newtheorem{rem}[thm]{Remark}
\theoremstyle{remark}

\theoremstyle{remark}

\theoremstyle{remark}
\newtheorem{question}[thm]{Question}
\theoremstyle{definition}
\newtheorem{example}[thm]{Example}

\theoremstyle{definition}

\theoremstyle{plain}

\theoremstyle{definition}

\theoremstyle{remark}

\theoremstyle{remark}

\theoremstyle{definition}

\theoremstyle{remark}
  \newtheorem*{acknowledgement*}{Acknowledgement}

\newcommand{\R}{{\mathbb R}}

\newcommand{\C}{{\mathbb C}}
\newcommand{\N}{{\mathbb N}}
\newcommand{\Z}{{\mathbb Z}}
\newcommand{\Q}{{\mathbb Q}}

\newcommand{\id}{\mathrm{id}}

\DeclareMathOperator{\Tr}{Tr}

\def\freeprod{\font\bigsymbolsfont=cmsy10 scaled \magstep3
 \setbox0=\hbox{\bigsymbolsfont\char'003 }\mathord{\lower1pt\box0}}\relax\ignorespaces

\newcommand{\Hawaii}{Hawai\kern.05em`\kern.05em\relax i}

\usepackage{amssymb} \usepackage{amscd} \usepackage{hyperref}  \usepackage{tikz} \usetikzlibrary{matrix,arrows} 
\usepackage{soul,xcolor} \usepackage{color} \usepackage{multirow}

\makeatother

\begin{document}
	
\setstcolor{red}

\title{Classification of Uniform Roe algebras of locally finite groups} 

\author{Kang Li$^{1}$}
\address{Mathematisches Institut der WWU M\"unster,
\newline Einsteinstrasse 62, 48149 M\"unster, Deutschland}
\email{lik@uni-muenster.de}

\author{Hung-Chang Liao$^{2}$}
\address{Mathematisches Institut der WWU M\"unster,
\newline Einsteinstrasse 62, 48149 M\"unster, Deutschland}
\email{liao@uni-muenster.de}

\thanks{{$^{1}$} Supported by the Danish Council for Independent Research (DFF-5051-00037) and is partially supported by the DFG (SFB 878).}

\thanks{{$^{2}$} Supported by Deutsche Forschungsgemeinschaft (SFB 878, Groups, Geometry and Actions).}

\date{\today}
\maketitle

\begin{abstract}  
	We study the uniform Roe algebras associated to locally finite groups. We show that for two countable locally finite groups $\Gamma$ and $\Lambda$, the associated uniform Roe algebras $C^*_u(\Gamma)$ and $C^*_u(\Lambda)$ are $*$-isomorphic if and only if their $K_0$ groups are isomorphic as ordered abelian groups with units. This can be seen as a non-separable non-simple analogue of the Glimm-Elliott classification of UHF algebras. To the best of our knowledge, this is the first classification result for a class of non-separable unital $C^*$-algebras. Along the way we also obtain a rigidity result: two countable locally finite groups are bijectively coarsely equivalent if and only if the associated uniform Roe algebras are $*$-isomorphic. 
	
	Finally, we give a summary of $C^*$-algebraic characterizations for (not necessarily countable) locally finite discrete groups in terms of their uniform Roe algebras. In particular, we show that a discrete group $\Gamma$ is locally finite if and only if the associated uniform Roe algebra $\ell^\infty(\Gamma)\rtimes_r \Gamma$ is locally finite-dimensional.
\end{abstract}

\section{Introduction}

Given a countable discrete group $\Gamma$, one can always equip $\Gamma$ with a proper left-invariant metric $d$ and such metric $d$ is unique up to bijective coarse equivalence (see \cite[Lemma 2.1]{Tu01}). In this way one obtains a canonical metric space structure $(\Gamma, d)$ on the group. To the metric space $(\Gamma,d)$ one can associate a $C^*$-algebra $C^*_u(\Gamma)$, called the \emph{uniform Roe algebra}, which encodes many large-scale properties of the group $\Gamma$. For instance, $C^*_u(\Gamma)$ is nuclear if and only if the metric space $(\Gamma, d)$ has Yu's property A if and only if the group $\Gamma$ is exact (see \cite{Oza00, MR1739727}). For amenability, it is shown in \cite{MR2873171} that $C^*_u(\Gamma)$ is a properly infinite $C^*$-algebra if and only if $\Gamma$ is non-amenable (more characterizations of proper infiniteness of uniform Roe algebras will be contained in \cite{ALLW17}; see also \cite[Theorem 4.2]{LW}). Another result along the same line is the work of Kellerhals, Monod, and R\o rdam on supramenability \cite{KMR13}. They showed that $\Gamma$ is supramenable if and only if $C^*_u(\Gamma)$ contains no properly infinite projection. 

More recently, Scarparo proved in \cite{Sca17} that the group $\Gamma$ is locally finite if and only if the uniform Roe algebra $C^*_u(\Gamma)$ is a finite $C^*$-algebra. By definition, a group is called \emph{locally finite} if all of its finitely generated subgroups are finite. For a countable discrete group $\Gamma$, a result of Smith asserts that local finiteness corresponds precisely to having asymptotic dimension zero (see \cite[Theorem 2]{Smi06}). In this case, it follows from the work of Winter and Zacharias on nuclear dimension of $C^*$-algebras \cite{WZ10} that its uniform Roe algebra $C^*_u(\Gamma)$ is \emph{locally finite-dimensional} (or \emph{local AF}), in the sense that given any finite subset $\mathcal{F}$ of $C^*_u(\Gamma)$ there is a finite-dimensional $C^*$-subalgebra of $C^*_u(\Gamma)$ which almost contains $\mathcal{F}$. In view of the classical Glimm-Elliott classification of UHF and (local) AF algebra (\cite{Gli60, Ell76}) using $K$-theory, one may ask whether $K$-theory carries useful information of the $C^*$-algebra $C^*_u(\Gamma)$ when it is locally finite-dimensional. In this paper we explicitly compute the $K$-theory of uniform Roe algebras $C^*_u(\Gamma)$ of countable locally finite groups, and show that in fact in this case $K$-theory is a complete invariant. To the best of our knowledge, this is the first classification result for a class of non-separable unital $C^*$-algebras. Moreover, the $K$-theory completely encodes geometric information of the underlying metric spaces. The following theorem is the main result of this paper.

\begin{thm}[Theorem \ref{thm:main}] \label{thm:main_intro}
	Let $\Gamma$ and $\Lambda$ be countable locally finite groups with proper left-invariant metrics $d_\Gamma$ and $d_\Lambda$, respectively. Then the following are equivalent:
	\begin{enumerate}
		\item $(\Gamma, d_\Gamma)$ and $(\Lambda, d_\Lambda)$ are bijectively coarsely equivalent;
		\item there is a $*$-isomorphism $\varphi:C^*_u(\Gamma)\to C^*_u(\Lambda)$ such that $\varphi(\ell^\infty(\Gamma)) = \ell^\infty(\Lambda)$; 
		\item $C^*_u(\Gamma)$ and $C^*_u(\Lambda)$ are $*$-isomorphic;
		\item $(K_0(C^*_u(\Gamma)), K_0(C^*_u(\Gamma))^+, [1]_0 ) \cong (K_0(C^*_u(\Lambda)), K_0(C^*_u(\Lambda))^+, [1]_0 )$;
		\item $( K_0(C^*_u(\Gamma)), [1]_0  ) \cong ( K_0(C^*_u(\Lambda)), [1]_0   )$. 
	\end{enumerate}
\end{thm}
It is already known that (1) and (2) are equivalent for all countable discrete groups (see Remark~\ref{(1)=(2)}). Note that the implication $(3)\Longrightarrow (1)$ reveals a rigidity phenomenon. More precisely, given a countable locally finite group $\Gamma$ and any countable discrete group $\Lambda$, if we know $C^*_u(\Gamma)$ and $C^*_u(\Lambda)$ are *-isomorphic, then $\Lambda$ must be bijectively coarsely equivalent to $\Gamma$ (in particular, $\Lambda$ must be locally finite). This type of rigidity result has been studied by \v{S}parkula and Willett in \cite{SW13}, where they showed, among other things, that the implication $(3)\Longrightarrow (1)$ holds for the class of non-amenable exact countable groups (see Remark~\ref{rigidity} for more details). Theorem~\ref{thm:main_intro} provides the same kind of rigidity result for a small class of amenable groups.\\

In addition to the study of $K$-theory, we are also interested in $C^*$-algebraic characterizations of locally finite groups in terms of their uniform Roe algebras. It is well-known that for a countable group $\Gamma$, the uniform Roe algebra $C^*_u(\Gamma)$ is $*$-isomorphic to the reduced crossed product $C^*$-algebra $\ell^\infty(\Gamma)\rtimes_r \Gamma$, where $\Gamma$ acts on $\ell^\infty(\Gamma)$ by right-translations (see e.g. \cite[Proposition 5.1.3]{BO08}).
Motivated by the above-mentioned result of Winter-Zacharias, we provide the following proposition:

\begin{prop}[Proposition \ref{local AF}] \label{prop:local_AF_intro}
	Let $\Gamma$ be a locally finite discrete (not necessarily countable) group. Then the reduced crossed product $\ell^\infty(\Gamma)\rtimes_r\Gamma$ is locally finite-dimensional.
\end{prop}

Combining this proposition with the result of Scarparo on finiteness of uniform Roe algebras \cite{Sca17}, we obtain the following list of characterizations:

\begin{cor}[Corollary \ref{cor:structure}] \label{thm:structure_intro}
	Let $\Gamma$ be a discrete group. Then the following are equivalent:
	\begin{enumerate}
		
		\item $\ell^\infty(\Gamma)\rtimes_r\Gamma$ is locally finite-dimensional;
		\item $\ell^\infty(\Gamma)\rtimes_r\Gamma$ is quasidiagonal;
		\item $\ell^\infty(\Gamma)\rtimes_r\Gamma$ has stable rank one;
		\item $\ell^\infty(\Gamma)\rtimes_r\Gamma$ has cancellation;
		\item $\ell^\infty(\Gamma)\rtimes_r\Gamma$ is stably finite;
		\item $\ell^\infty(\Gamma)\rtimes_r\Gamma$ is finite;
		\item $\Gamma$ is locally finite.
	\end{enumerate}
\end{cor}

We should mention that the analogous result stated in Corollary~\ref{thm:structure_intro} does not hold for general metric spaces (see Remark \ref{rem:structure_space}). 

Let us briefly describe how this piece is organized. In Section 2 we recall some basic definitions and results from coarse geometry, and establish our notations. In particular, we discuss the notion of large-scale connected components, which plays a crucial role in our study later on. In Section 3 we focus on the class of locally finite countable groups and explore the connection between their large-scale geometric properties and group-theoretic properties. In particular, the bijective coarse equivalence class of a countable locally finite group $\Gamma$ is completely determined by the cardinality of Sylow subgroups of $\Gamma$. In Section 4 we compute the $K$-theory of uniform Roe algebras built from countable locally finite groups and prove our main theorem. This is done by making use of the inductive limit decomposition of the uniform Roe algebra and a detailed study of the connecting maps. The last section is devoted to various equivalent $C^*$-algebraic properties of uniform Roe algebras of (not necessarily countable) locally finite discrete groups.


\section{Preliminaries}

In this section we review some basic definitions and constructions from coarse geometry and the associated $C^*$-algebras.  Let $(X,d_X)$ and $(Y,d_Y)$ be metric spaces. A map $f:X\to Y$ is called \emph{bornologous} (or \emph{uniformly expansive}) if for every $R>0$ there exists $S>0$ such that $d_X(x,x') \leq R$ implies $d_Y(f(x),f(x')) \leq S$. Two maps $f,f':X\to Y$ are said to be \emph{close} if the function $x\mapsto d_Y(f(x),f'(x))$ is bounded on $X$. A bornologous map $f:X\to Y$ is called a \emph{coarse equivalence} if there exists a bornologous map $g:Y\to X$ such that $f\circ g$ is close to $\id_Y$ and $g\circ f$ is close to $\id_X$.

\begin{defn}
	Let $(X,d_X)$ and $(Y,d_Y)$ be metric spaces. We say a map $f:X\to Y$ is a \emph{bijective coarse equivalence} if $f$ is both a coarse equivalence and a bijection. In this case we say $X$ and $Y$ are \emph{bijectively coarsely equivalent}.
\end{defn}

A metric space $(X,d)$ is said to have \emph{bounded geometry} if for every $r>0$, the function $x\mapsto |B(x,r)|$ is bounded on $X$, where $B(x,r) = \{ x'\in X: d(x,x') \leq r \}$ denotes the ball with center $x$ and radius $r$. A metric space with bounded geometry is necessarily discrete and countable.

\begin{defn}
	Let $(X,d)$ be a metric space with bounded geometry. An operator $T$ in $B(\ell^2(X))$ is said to have \emph{finite propagation} if there exists $R>0$ such that $\langle T\delta_{x'}, \delta_x\rangle = 0$ whenever $d(x,x') > R$.
	
	Let $\C_u[X]$ be the $*$-subalgebra of all operators in $B(\ell^2(X))$ with finite propagation. The norm completion of $\C_u[X]$ in $B(\ell^2(X))$ is called the \emph{uniform Roe algebra} of $(X,d)$, which is denoted $C^*_u(X)$.
\end{defn}

Note that for a metric space $(X,d)$, the uniform Roe algebra $C^*_u(X)$  contains all compact operators on the Hilbert space $\ell^2(X)$. Moreover, it contains $\ell^\infty(X)$ as diagonal matrices. Therefore, $C^*_u(X)$ is non-simple and non-separable for every infinite metric space $X$.

It is well-known that the uniform Roe algebra is invariant under bijective coarse equivalence. We state the following proposition and give a proof mainly for the reader's convenience.

\begin{prop} \cite[Proof of Theorem 4]{BNW07} \label{prop:bce}
	Let $(X,d_X)$ and $(Y,d_Y)$ be metric spaces with bounded geometry. If $X$ and $Y$ are bijectively coarsely equivalent, then there is a $*$-isomorphism $\varphi:C^*_u(X)\to C^*_u(Y)$ such that $\varphi(\ell^\infty(X)) = \ell^\infty(Y)$. In particular, $C^*_u(X)$ and $C^*_u(Y)$ are $*$-isomorphic.
\end{prop}
\begin{proof}
	Let $f:X\to Y$ be a bijective coarse equivalence. Then there is a unitary $u:\ell^2(X)\to \ell^2(Y)$ satisfying $u\delta_x = \delta_{f(x)}$ for all $x\in X$. Define
	$$
	\varphi:B(\ell^2(X))\to B(\ell^2(Y)),\;\;\;\;\;\; \varphi(T) = uTu^*.
	$$
	We show that $\varphi$ maps $C^*_u(X)$ into $C^*_u(Y)$. First note that for all $x_1,x_2$ in $X$, we have
	\begin{align*}
	\langle \varphi(T)\delta_{f(x_2)}, \delta_{f(x_1)} \rangle = \langle uTu^*(\delta_{f(x_2)}), \delta_{f(x_1)} \rangle = \langle T\delta_{x_2}, \delta_{x_1} \rangle.	
	\end{align*}
	Let $T$ be an operator in $B(\ell^2(X))$ of finite propagation, say $\langle T\delta_{x'}, \delta_x \rangle = 0$ whenever $d(x',x) > R$. Since $f$ is bornologous, there exists $S>0$ such that
	$d_Y(f(x_1),f(x_2))<S$ whenever $d(x_1,x_2) <R$. Then it is clear from the above calculation that $\varphi(T)$ has finite propagation. Since $\varphi$ is an isometry, it actually maps the entire $C^*_u(X)$ into $C^*_u(Y)$.
	
	By the same reasoning the map $\psi:T\mapsto u^*Tu$ maps $C^*_u(Y)$ into $C^*_u(X)$. It is clear $\psi$ is the inverse of $\varphi$, and hence $\varphi$ is an isomorphism. Moreover, it is readily seen from the definition that $\varphi$ maps $\ell^\infty(X)$ onto $\ell^\infty(Y)$. 
\end{proof}

Below we recall the notion of large-scale connectedness. This would become important later as it allows us to decompose the uniform Roe algebra of locally finite groups into an inductive limit, which in turn facilitates the computation of the $K$-theory.

\begin{defn} Let $(X,d)$ be a metric space and $R > 0$. Two elements $x$ and $y$ in $X$ are called $\emph{R-connected}$ if there exists a finite sequence $x_0,x_1,...,x_n$ in $X$ such that $x=x_0$, $y = x_n$, and $d(x_i,x_{i+1}) \leq R$ for all $i=0,1,...,n-1$. 
\end{defn}

This is easily seen to be an equivalence relation on $X$, and the equivalence classes are called the \emph{R-connected components} of $X$. For convenience we also talk about the \emph{$0$-connected components} of $X$, which are nothing but the points in $X$. The next lemma says that a bijective coarse equivalence behaves well in terms of large-scale connected components.

\begin{lem} \label{lem:component-bijection}
	Let $f:(X,d_X)\to (Y,d_Y)$ be a bijective coarse equivalence between metric spaces. Let $R> 0$ and $X = \bigsqcup_{i\in I_R} X_{R,i}$ be the decomposition of $X$ into $R$-connected components. Then there exists $S>0$ such that each $f(X_{R,i})$ is contained in some $S$-connected component of $Y$.
	
	Moreover, in this case each $S$-connected component $Y_{S,j}$ of $Y$ is a disjoint union $\bigsqcup_{i\in I_{R,j}} f(X_{R,i})$ for some subset $I_{R,j}$ of $I_R$.
\end{lem}
\begin{proof}
	The first statement follows from the definition and the fact that $f$ is bornologous. For the second statement, let $y$ be a point in $Y_{S,j}$. Since $f$ is surjective, there is a $R$-connected component $X_{R,i}$ of $X$ such that $f(X_{R,i})$ contains $y$. By the first statement $f(X_{R,i})$ is entirely contained in $Y_{S,j}$. Note that the union is disjoint because $f$ is injective.
\end{proof}

In this paper we are mainly interested in spaces (in fact, groups) with zero dimension in the large scale sense. The following definition is slightly non-standard, but is readily seen to be equivalent to the usual one.

\begin{defn}
	Let $(X,d)$ be a metric space. We say $X$ has \emph{asymptotic dimension zero} if for every $R>0$, there is a uniform bound on the diameters of the $R$-connected components of $X$.
\end{defn}

Note that when $(X,d)$ has bounded geometry, having asymptotic dimension zero is equivalent to having a uniform bound on the cardinalities of the $R$-connected components of $X$.



\section{Locally finite groups and supernatural numbers}

Starting this section we will mainly focus on locally finite groups. As we will see, this is precisely the class of groups which have asymptotic dimension zero.

\begin{defn}
	A group $G$ is called \emph{locally finite} if every finitely generated subgroup of $G$ is finite.
\end{defn}

The following (easy) characterization of local finiteness will become handy later.

\begin{lem} \label{lem:increase-sequence}
	Let $\Gamma$ be a countable group. Then $\Gamma$ is locally finite if and only if there exists an increasing sequence 
	$$
	\{e\} =: \Gamma_0 \subseteq \Gamma_1 \subseteq \Gamma_2 \subseteq \cdots \subseteq \Gamma
	$$
	of finite subgroups of $\Gamma$ such that $\Gamma = \bigcup_{n=0}^\infty \Gamma_n$.
\end{lem}

\begin{example} Here we list some examples of locally finite groups.
	\begin{enumerate}
		\item Every finite group is clearly locally finite.
		\item Every infinite direct sum of finite groups, such as $\bigoplus_{\mathbb{N}}\Z/2\Z$, is locally finite.
		\item $S_\infty$, the direct limit of all finite permutation groups, is locally finite.
		\item $\Q/\Z$ is locally finite.
	\end{enumerate}
\end{example}

Recall that a metric $d$ on a discrete group $\Gamma$ is called \emph{proper} if any closed bounded subset of $\Gamma$ is finite, and $d$ is \emph{left-invariant} if $d(s,t)$ = $d(gs,gt)$ for all $g,s,t$ in $\Gamma$. Every countable discrete group $\Gamma$ can be equipped with a proper left-invariant metric $d$ in the following way: let
$$
\{e\} =: F_0 \subseteq F_1 \subseteq F_2 \subseteq \cdots \subseteq \Gamma
$$
be an increasing sequence of finite symmetric subsets of $\Gamma$ such that $\Gamma = \bigcup_{n=0}^\infty F_n$ and $F_nF_m\subseteq F_{n+m}$ for every $n, m\in \N$. Define $d:\Gamma\times \Gamma\to \R$ by
$$
d(g,h) := \min\{ n\in \N: g^{-1}h \in F_n  \}\;\;\;\;\;\; (g,h\in \Gamma).
$$
It is routine to check that $d$ is a proper left-invariant metric on $\Gamma$.

\begin{prop} \cite[Lemma 2.1]{Tu01}
	Let $\Gamma$ be a countable discrete group with proper left-invariant metrics $d_1$ and $d_2$. Then the identity map
	$$
	\id:(\Gamma, d_1) \to (\Gamma, d_2)
	$$
	is a (bijective) coarse equivalence.
	
	As a consequence, there is a unique proper left-invariant metric on $\Gamma$ up to bijective coarse equivalence.
\end{prop}

Now let $\Gamma$ be a countable locally finite group. By Lemma \ref{lem:increase-sequence} there is an increasing sequence of finite subgroups
$$
\{e\} =: \Gamma_0\subseteq \Gamma_1 \subseteq \Gamma_2\subseteq \cdots \subseteq \Gamma
$$
such that $\Gamma = \bigcup_{n=0}^\infty \Gamma_n$. As we saw in the previous paragraph, this sequence gives rise to a proper left-invariant metric $d$ on $\Gamma$. In this case we have a very simple description of the $n$-connected components.

\begin{lem} \label{lem:group-component}
	Let $(\Gamma, \{\Gamma_n\}_{n=0}^\infty,\; d)$ be as above. Then the $n$-connected components ($n=0,1,2,...$) are exactly the left cosets of $\Gamma_n$.
\end{lem}
\begin{proof}
	Observe that $d(g,h) \leq n$ if and only if $g$ and $h$ are in the same left coset of $\Gamma_n$. The result then follows.
\end{proof}

With this observation it is clear that locally finite groups have asymptotic dimension zero. In fact the converse is also true, as shown in the following theorem of Smith:

\begin{thm} \cite[Theorem 2]{Smi06} \label{thm:char-loc-fin}
	Let $\Gamma$ be a countable group with a proper left-invariant metric $d$. Then the following are equivalent:
	\begin{enumerate}
		\item $\Gamma$ is locally finite;
		\item $\Gamma$ has asymptotic dimension zero.
	\end{enumerate}
\end{thm}

Our next goal is to undertand a result of Protasov (see Theorem \ref{thm:supernatural}). Roughly speaking, the result says that the bijective corase equivalence class of a countable locally finite group $\Gamma$ is determined by the cardinality of the Sylow subgroups of $\Gamma$. This result was originally stated in terms of what is called ``ball's structure''. The aim here is to rephrase it in our language, and for the reader's convenience we also give a short proof of the theorem.

Let $\Gamma$ be a countable locally finite group. Let $\{p_1,p_2,...\}$ be the set of all prime numbers listed in increasing order. For each $j$ define
$$
n_j := \sup \{ m\in \N: p_j^m \text{ divides $|F|$ for some finite subgroup $F$ of $\Gamma$ } \}.
$$
Then the sequence $\{n_j\}_{j=1}^\infty$ is a so-called \emph{supernatural number}, which we denote by $s(\Gamma)$.  We usually think of a supernatural number $\{n_j\}_{j=1}^\infty$ as a formal product $p_1^{n_1}p_2^{n_2}\cdots$. Therefore for each prime number $p_j$ in the list, we say $p_j^m$ \emph{divides} $s(\Gamma) = \{n_j\}_{j=1}^\infty$ if $m\leq n_j$. Two supernatural numbers are \emph{equal} if they are equal as sequences.

\begin{example}
	\begin{enumerate} 
		\item For any prime number $p$, we have $s(\bigoplus_{\N}\Z/p\Z  ) = p^\infty$.
		\item Since the cardinality of the permutation group $S_n$ of $n$ elements is equal to $(n!)$, we see that $s(S_\infty) = 2^\infty3^\infty5^\infty\cdots$.
		\item Given any prime number $p$ and natural number $k$, the subgroup of $\Q/\Z$ generated by $\frac{1}{p^k}$ is isomorphic to $\Z/p^k\Z$. Therefore $s(\Q/\Z) = 2^\infty3^\infty5^\infty\cdots$.
	\end{enumerate}
\end{example}

\begin{prop}
	Let $\Gamma$ and $\Lambda$ be two countable locally finite groups equipped with proper left-invariant metrics $d_\Gamma$ and $d_\Lambda$, respectively. Suppose there is a bijective coarse equivalence $f:\Gamma\to \Lambda$. Then for every finite subgroup $F$ of $\Gamma$ there is a finite subgroup $E$ of $\Lambda$ such that $|F|$ divides $|E|$.
\end{prop}
\begin{proof}
	We may assume that the metric $d_\Gamma$ comes from an increasing sequence $\{\Gamma_n\}_{n=0}^\infty$ of finite subgroups such that $\bigcup_{n=0}^\infty \Gamma_n = \Gamma$, and similarly $d_\Lambda$ comes from an increasing sequence $\{\Lambda_n\}_{n=0}^\infty$ of finite subgroups of $\Lambda$. Given a finite subgroup $F$ of $\Gamma$, there is some $n\in \N$ such that $F$ is contained in $\Gamma_n$. Recall that the $n$-connected components of $\Gamma$ are exactly the left cosets of $\Gamma_n$. By Lemma \ref{lem:component-bijection}, $f(\Gamma_n)$ is contained in an $m$-connected component of $\Lambda$ for some $m$, which is a left coset of $\Lambda_m$. Moreover this left coset of $\Lambda_m$ is a disjoint union of the images of some left cosets of $\Gamma_n$. Therefore $|\Gamma_n|$ (and hence $|F|$) divides $|\Lambda_m|$.
\end{proof}

\begin{cor}
	Let $\Gamma$ and $\Lambda$ be two countable locally finite groups equipped with proper left-invariant metrics $d_\Gamma$ and $d_\Lambda$, respectively. If $(\Gamma,d_\Gamma)$ and $(\Lambda, d_\Lambda)$ are bijectively coarsely equivalent, then $s(\Gamma) = s(\Lambda)$.
\end{cor}

We have seen that the supernatural number $s(\Gamma)$ is an invariant under bijective coarse equivalence. The next theorem asserts that it is actually a complete invariant. For those who are familiar with the theory of UHF algebras, this result strongly resembles Glimm's classification theorem of UHF algebras.

\begin{thm} \cite[Theorem 5]{Pro02} \label{thm:supernatural}
	Let $\Gamma$ and $\Lambda$ be two countable locally finite groups with proper left-invariant metrics $d_\Gamma$ and $d_\Lambda$, respectively. Then the following are equivalent.
	\begin{enumerate}
		\item $(\Gamma, d_\Gamma)$ and $(\Lambda, d_\Lambda)$ are bijectively coarsely equivalent.
		\item $\Gamma$ and $\Lambda$ have the same supernatural number, i.e., $s(\Gamma) = s(\Lambda)$.
	\end{enumerate}

		In particular, there are uncountably many bijective coarse equivalence classes of countable locally fintie groups.
\end{thm}
\begin{proof}
	It remains to prove that $(2)$ implies $(1)$. As before we may assume that the metric $d_\Gamma$ comes from an increasing sequence
	$$
	\{e_\Gamma\} = \Gamma_0 \subseteq \Gamma_1 \subseteq \Gamma_2\subseteq \cdots \subseteq \Gamma
	$$
	of subgroups of $\Gamma$ and similarly $d_\Lambda$ comes from
	$$
	\{e_\Lambda\} = \Lambda_0 \subseteq \Lambda_1 \subseteq \Lambda_2\subseteq \cdots \subseteq \Lambda.
	$$
	For convenience we may also assume without losing generality that
	\begin{enumerate}
		\item [(a)]$|\Gamma_k|$ divides $|\Lambda_k|$ and $|\Lambda_k|$ divides $|\Gamma_{k+1}|$ for all $k\in \N$, and
		\item [(b)]$|\Gamma_k| \neq |\Lambda_k|$ and $|\Lambda_k| \neq |\Gamma_{k+1}|$ for all $k\in \N$.
	\end{enumerate}  
	
	We will construct, inductively, maps $\varphi_k:\Gamma_k\to \Lambda_k$ and $\psi_k:\Lambda_k\to \Gamma_{k+1}$ for $k=0,1,2,...$ such that
	\begin{enumerate}
		\item $\varphi_k \equiv \varphi_{k+1}$ on $\Gamma_k$ for each $k=0,1,2,...$;
		\item $\psi_k \equiv \psi_{k+1}$ on $\Lambda_k$ for each $k=0,1,2,...$;
		\item $\psi_k\circ \varphi_k(x) = x$ for all $k\in \N$ and all $x\in \Gamma_k$;
		\item $\varphi_{k+1}\circ \psi_k (y) = y$ for all $k\in \N$ and all $y\in \Lambda_k$;
		\item Given any $\ell\leq k$ and $x_1,x_2$ in $\Gamma_k$ such that $x_1$ and $x_2$ belong to the same left coset of $\Gamma_{\ell}$ in $\Gamma_k$, the images $\varphi_k(x_1)$ and $\varphi_k(x_2)$ belong to the same left coset of $\Lambda_{\ell}$ in $\Lambda_k$.
		\item Given any $\ell\leq k$ and $y_1,y_2$ in $\Lambda_k$ such that $y_1$ and $y_2$ belong to the same left coset of $\Lambda_{\ell}$ in $\Lambda_k$, the images $\psi_k(y_1)$ and $\psi_k(y_2)$ belong to the same left coset of $\Gamma_{\ell+1}$ in $\Gamma_{k+1}$.
	\end{enumerate}
	Once these maps are constructed, we define
	$$
	\varphi:\Gamma\to \Lambda,\;\;\;\;\;\; \varphi(x) := \varphi_k(x)\;\;\;\;\;\; (x\in \Gamma_k)
	$$
	and similarly
	$$
	\psi:\Lambda\to \Gamma,\;\;\;\;\;\; \psi(y) := \psi_k(y)\;\;\;\;\;\; (y\in \Lambda_k).
	$$
	By (1) and (2) these two maps are well-defined. The map $\varphi$ is bijective because of (3) and (4) (and $\psi$ is precisely the inverse of $\varphi$), and from (5) and (6) one easily checks that $\varphi$ is a coarse equivalence.
	
	It remains to construct the maps $\{\varphi_k\}_{k=0}^\infty$ and $\{\psi_k\}_{k=0}^\infty$. Define $\varphi_0:\Gamma_0\to \Lambda_0$ by $\varphi_0(e_\Gamma) = e_\Lambda$ and $\psi_0:\Lambda_0\to \Gamma_1$ by $\psi_0(e_\Lambda) = e_\Gamma$. Now suppose we have constructed the maps $\varphi_0,\varphi_1,...,\varphi_{k-1}$ and $\psi_0,\psi_1,...,\psi_{k-1}$. We need to define $\varphi_k$ and $\psi_k$. Since $|\Lambda_{k-1}|$ divides $|\Gamma_k|$, we have $|\Gamma_k| = m|\Lambda_{k-1}|$ for some $m\in \N$. Choose $m$ distinct left cosets of $\Lambda_{k-1}$ in $\Lambda_k$, say $s_0\Lambda_{k-1}, s_1\Lambda_{k-1},...,s_{m-1}\Lambda_{k-1}$ ($s_0,s_1,...,s_{m-1}\in \Lambda_k$), with one of them being the subgroup $\Lambda_{k-1}$ itself. Then there exists a map
	$$
	\varphi_k:\Gamma_k\to \Lambda_k
	$$
	such that 
	\begin{itemize}
		\item $\varphi_k$ agrees with $\varphi_{k-1}$ on $\Gamma_{k-1}$, and
		\item $\varphi_k$ maps $\Gamma_k$ bijectively onto the union $\bigsqcup_{j=0}^{m-1} s_j\Lambda_{k-1}$ in the way that condition (5) holds (this is possible because of the divisibility assumption among the subgroups $\Gamma_\ell$ and $\Lambda_\ell$).
	\end{itemize}
	This completes the construction of $\varphi_k$. Now $\psi_k$ is defined in a completely analogous way.
\end{proof}

\begin{rem}\label{BZ}
	In contrast to Theorem \ref{thm:supernatural}, Banakh and Zarichnyi showed in \cite[Corollary~8]{BZ11} that all countably infinite locally finite groups are coarsely equivalent. 
\end{rem}


\section{The K-theory and classification}

In this section we prove our main theorem (Theorem \ref{thm:main_intro}). We begin by computing the $K_0$-group of an infinite product of matrix algebras. This result is well-known to experts, and a proof is included only for the reader's convenience.

\begin{prop} \label{prop:K0-propduct} Suppose $\{s_i\}_{i=1}^\infty$ is a bounded sequence of positive integers. Let $A = \prod_{i=1}^\infty M_{s_i}(\mathbb{C})$. Then
	$$
	(K_0( A ), K_0( A )^+, [1_A]_0  ) \cong (\ell^\infty(\N,\Z), \ell^\infty(\N,\Z)^+, \{s_i\}_{i=1}^\infty),
	$$
	where $\ell^\infty(\N,\Z)^+$ is the set of all positive sequences in $\ell^\infty(\N,\Z)$.
\end{prop}
\begin{proof}
	For any matrix algebra $M_n(\C)$ we write $\Tr$ for the non-normalized trace on $M_n(\C)$. Let $P_n(A)$ be the set of projections in $M_n(A)$ and write $P_\infty(A) = \bigcup_{n=1}^\infty P_n(A)$. Consider the map
	$$
	P_\infty(A)\to \Z^\N,\;\;\;\;\;\;	\mathbf{p} \mapsto \{ \Tr(p_i) \}_{i=1}^\infty,
	$$
	where $\mathbf{p} = (p_1,p_2,...)$ is a projection in $M_n(A) \cong \prod_{i=1}^\infty M_{ns_i}(\C)$. This induces a well-defined group homomorphism
	$$
	\Tr_*:K_0(A)\to \Z^\N
	$$
	which satisfies the formula
	$$
	\Tr_*([\mathbf{p}]_0) = \{ \Tr(p_i) \}_{i=1}^\infty.
	$$
	Since the sequence $\{s_i\}_{i=1}^\infty$ is bounded, the image of $\Tr_*$ is contained in the subgroup $\ell^\infty(\N,\Z)$ of $\Z^\N$. We would like to show that $\Tr_*$ is an isomorphism onto $\ell^\infty(\N,\Z)$.
	
	First of all, if $\Tr_*( [\mathbf{p}]_0 - [\mathbf{q}]_0  ) = 0$ for some projections $\mathbf{p},\mathbf{q}$ in $M_n(A)$, then $\Tr(p_i) = \Tr(q_i)$ for all $i=1,2,...$. Therefore $p_i$ is Murray-von Neumann equivalent to $q_i$ in $M_{ns_i}(\C)$ for each $i=1,2,...$. It follows that $\mathbf{p}$ and $\mathbf{q}$ are Murray-von Neumann equivalent in $M_n(A)$ and hence $\Tr_*$ is injective.
	
	As for surjectivity, it suffices to show that the image of $\Tr_*$ contains $\ell^\infty(\N,\Z)^+$. Suppose we are given a sequence $\{r_i\}_{i=1}^\infty$ in $\ell^\infty(\N,\Z)^+$ which is bounded by $N\in \N$. Let $\mathbf{p} = \{p_i\}$ be a projection in $M_N(A)$ such that each $p_i$ has rank $r_i$. Then clearly we have $\Tr_*([\mathbf{p}]_0) = \{r_i\}_{i=1}^\infty$. Hence $\Tr_*$ is surjective.
	
	Note that the previous paragraph also shows that $\Tr_*$ maps the positive cone $K_0(A)^+$ onto $\ell^\infty(\N,\Z)^+$. Finally it is clear from the definition that $\Tr_*$ maps the class $[1_A]_0$ to the sequence $\{s_i\}_{i=1}^\infty$.
\end{proof}

Let us return briefly to general metric spaces. For the following discussion we assume the metric space $(X,d)$ has bounded geometry and asymptotic dimension zero. For each $R>0$ let $X = \bigsqcup_{i\in I_R}X_{R,i}$ be the decomposition of $X$ into $R$-connected components.  Suppose $0<R<R'$. Since each $R'$-connected component is a disjoint union of $R$-connected components, there is a canonical embedding
$$
\varphi_{R,R'}:\prod_{i\in I_R} B(\ell^2(X_{R,i})) \to \prod_{i\in I_{R'}} B(\ell^2(X_{R',i})).
$$
as block-diagonal matrices. Moreover, since $X$ has asymptotic dimension zero,  for each $R>0$ there is a natural inclusion
$$
\mu_R: \prod_{i\in I_R} B(\ell^2(X_{R,i})) \to C^*_u(X)
$$
(again as block-diagonal matrices) such that the following diagram
\begin{center}
	\begin{tikzpicture} [node distance = 5cm, auto]
	\node (1) {$\prod_{i\in I_R} B(\ell^2(X_{R,i})) $};
	\node (2) [right of = 1] {$\prod_{i\in I_{R'}} B(\ell^2(X_{R',i}))$};
	\node (3) [right of = 1, node distance = 2.5cm, below of = 1] {$C^*_u(X)$};
	\draw [->] (1) to node {$\varphi_{R,R'}$} (2);
	\draw [->] (2) to node {$\mu_{R'}$} (3);
	\draw [->] [swap] (1) to node {$\mu_{R}$} (3);
	\end{tikzpicture}
\end{center}
commutes for all $R' > R$.

Recall that we also talk about the $0$-connected components, i.e., the points of $X$. The corresponding algebra is then $\prod_{i\in I_0} B(\ell^2(X_{0,i})) \cong \prod_{x\in X} \mathbb{C}$. There are also natural embeddings of $\prod_{x\in X}\mathbb{C}$ into $\prod_{i\in I_R} B(\ell^2(X_{R,i}))$ and $C^*_u(X)$ as before.

\begin{prop} \label{prop:zero-dim-ind-lim}
Let $(X,d)$ be a metric space with bounded geometry. If $(X,d)$ has asymptotic dimension zero, then	
$$
C^*_u(X) \cong \varinjlim \left( \prod_{i\in I_n} B(\ell^2(X_{n,i})), \varphi_n \right)\;\;\;\;\;\; (n=0,1,2,...),
$$
where $\varphi_n := \varphi_{n,n+1}$.
\end{prop}
\begin{proof}
	It suffices to check that the union $\bigcup_{n=0}^\infty \mu_n\left( \prod_{i\in I_n} B(\ell^2(X_{n,i})) \right)$ contains all operators in $C^*_u(X)$ of finite propagation. Suppose $T$ is an operator on $\ell^2(X)$ which has finite propagation. Then by definition there exists $S>0$ such that $\langle T\delta_y, \delta_x \rangle = 0$ whenever $d(x,y) > S$. This implies that $T$ belongs to (the image of) $\prod_{i\in I_n} B(\ell^2(X_{n,i}))$ for any positive integer $n$ larger than $S$.	
\end{proof}

In order to write down the connecting maps more concretely, below we show that every metric space with asymptotic dimension zero is bijectively coarsely equivalent to a subspace of $\Z_{\geq 0}$, which we equip with the standard metric (see Proposition \ref{prop:subset-of-N}). Since $\Z_{\geq 0}$ has a natural order, this allows us to express the connecting maps in very simple terms (see Proposition \ref{prop:algebra-ind-lim}).

\begin{lem}
	Let $(X,d)$ be a metric space with bounded geometry which has asymptotic dimension zero. Fix $x_0$ in $X$ and write $X_n$ for the $n$-connected component of $X$ which contains $x_0$. Then for every $n\in \Z_{\geq 0}$, there exists a (finite) subset $Y_n$ of $\Z_{\geq 0}$ and a bijection $f_n:X_n\to Y_n$ such that $f_n(x_0) = 0$ and $f_n$ maps any $k$-connected component of $X_n$ onto a $k$-connected component of $Y_n$ (here $k\leq n$).
	
	Moreover, given a bijection $f_n:X_n\to Y_n$ as above, one can always choose $Y_{n+1}$ and $f_{n+1}:X_{n+1}\to Y_{n+1}$ in a way that $f_{n+1} \equiv f_n$ on $X_n$ and $Y_n$ is an $n$-connected component of $Y_{n+1}$.
\end{lem}
\begin{proof}
	We prove the lemma by induction on $n$. For $n=0$, define $Y_0 := \{0\}$ and let $f_0$ be the only possible map. Now assume the statement is true for $n=0,1,2,...,m-1$. Let $X_m = \bigsqcup_{i=0}^N X_{m-1,i}$ be the decomposition of $X_m$ into $(m-1)$-connected components of $X_m$ (these are also $(m-1)$-connected components of $X$) Note that since $X$ has asymptotic dimension zero, $X_m$ is finite and hence admits a finite decomposition. For each $i=0,1,...,N$ choose an element $z_i$ in $X_{m-1,i}$. Without losing generality we may assume that $X_{m-1,1} = X_{m-1}$ and $z_1 = x_0$. By the induction hypothesis we can find a subset $Y_{m-1,1}$ of $\Z_{\geq 0}$ and a bijection $g_1:X_{m-1,1}\to Y_{m-1,1}$ such that $g_1(z_1) = 0$ and $g_1$ maps any $k$-connected component of $X_{m-1,1}$ onto a $k$-connected component of $Y_{m-1,1}$ (here $k\leq m-1$). 
	
	Let $y_1 := \max\{p\in \Z_{\geq 0}: p\in Y_{m-1,1} \}$. In other words, let $y_1$ be the ``right-most'' point in $Y_{m-1,1}$. Using the induction hypothesis again, we find a subset $Y_{m-1,2}$ of $\Z_{\geq 0} \setminus \{ 0,1,...,y_1+{m-1} \}$ and a bijection $g_2:X_{m-1,2}\to Y_{m-1,2}$ such that $g_2(z_2) = y_1+m$ and $g_2$ satisfies the condition of mapping smaller connected components to connected components. Let $y_2 := \max\{ p\in \Z_{\geq 0}: p\in Y_{m-1,2} \}$ and proceed to define $Y_{m-1,3},...,Y_{m-1,N}$ and $g_3,...,g_N$ in the same way.	Now take $Y_m := \bigcup_{i=1}^N Y_{m,i}$ and define $f_m:X_m\to Y_m$ by pasting the maps $g_1,g_2,...,g_N$. It remains to observe that $Y_m$ is $m$-connected by construction.
	
	The second statement follows from the construction.
\end{proof}

\begin{prop} \label{prop:subset-of-N}
	Let $(X,d)$ be a metric space with bounded geometry which has asymptotic dimension zero. Then $X$ is bijectively coarsely equivalent to a subset of $\Z_{\geq 0}$ with respect to the standard metric.
\end{prop}
\begin{proof}
	Fix $x_0$ in $X$. For each $n\in \Z_{\geq 0}$ write $X_n$ for the $n$-connected component of $X$ which contains $x_0$. Then we have an increasing sequence
	$$
	\{x_0\} = X_0 \subseteq X_1\subseteq \cdots \subseteq X
	$$
	such that $X = \bigcup_{n=0}^\infty X_n$. Using the previous lemma one inductively constructs a sequence $\{Y_n\}$ of subsets of $\Z_{\geq 0}$ and a sequence of bijections $\{f_n:X_n\to Y_n\}$ such that
	\begin{enumerate}
		\item $f_{n+1}\equiv f_n$ on $X_n$,
		\item $Y_n$ is an $n$-connected component of $Y_{n+1}$, and
		\item $f_n$ maps any $k$-connected component of $X_n$ onto a $k$-connected component of $Y_n$ ($k\leq n$).
	\end{enumerate}
	Define 
	$$
	Y := \bigcup_{n=0}^\infty Y_n
	$$
	and 
	$$
	f:X\to Y,\;\;\;\;\;\; f(x) := f_n(x)\;\;\;\;\;\; (x\in X_n).
	$$
	By condition (1) the map $f$ is well-defined. Since each $f_n$ is a bijection, so is $f$. It is easy to check, using condition (3) and the fact that $X$ (and $Y$) have asymptotic dimension zero, that $f$ is a coarse equivalence.
\end{proof}

\begin{prop} \label{prop:algebra-ind-lim}
	Let $\Gamma$ be a countable locally finite group and let $(\Gamma, \{\Gamma_n\}_{n=0}^\infty,\; d)$ be as in Lemma \ref{lem:group-component}. For each $n=0,1,2,...$, define $k_n = |\Gamma_n|$ and $r_n = k_{n+1}/k_n$. Then 
	$$
	C^*_u(\Gamma) \cong \varinjlim \left( \prod_{i=1}^\infty M_{k_n}(\C), \varphi_n   \right),
	$$
	where 
	$$
	\varphi_n(T_1,T_2,...) = \left( \mathrm{diag}(T_1,...,T_{r_n}), \mathrm{diag}(T_{r_n+1},...,T_{2r_n}),...  \right).
	$$
\end{prop}
\begin{proof}
	We saw in Theorem \ref{thm:char-loc-fin} that locally finite groups have asymptotic dimension zero. Now the result  follows from Proposition \ref{prop:subset-of-N} and Proposition \ref{prop:zero-dim-ind-lim}.
\end{proof}

\begin{example} \label{example:Z_2}
	We consider the case when $\Gamma = \bigoplus_{\N} \Z/2\Z$. We can take the squence of finite subgroups $\{\Gamma_n\}$ to be $\Gamma_n = \bigoplus_{i=1}^n \Z/2\Z$. Then $k_n = |\Gamma_n| = 2^n$ and $r_n = k_{n+1}/k_n = 2$ for $n=0,1,2,...$. Then according to the previous proposition,
	$$
	C^*_u(\bigoplus_{\N}\Z/2\Z ) \cong \varinjlim \left( \prod_{i=1}^\infty M_{2^n}(\C), \varphi_n   \right),
	$$
	where
	$$
	\varphi_n(T_1,T_2,T_3,T_4,...) = \left(  \begin{pmatrix}
	T_1 & 0 \\
	0 & T_2
	\end{pmatrix},\begin{pmatrix}
	T_3 & 0 \\
	0 & T_4
	\end{pmatrix},...  \right).
	$$
\end{example}

From Propositionss \ref{prop:algebra-ind-lim} and Proposition \ref{prop:K0-propduct} we see that the ordered $K_0$-group of $C^*_u(\Gamma)$ is a (sequential) inductive limit of $\ell^\infty(\N,\Z)$. It is not hard to see that each connecting map $\varphi_n$ ($n=0,1,2,...$) induces the following map at the level of $K_0$-groups:
$$
\alpha_n :\ell^\infty(\N,\Z)\to \ell^\infty(\N,\Z),
$$
$$
\alpha_n( (m_1,m_2,...) ) =  (m_1+\cdots + m_{r_n},\; m_{r_n+1}+ \cdots + m_{2r_n},\;...).
$$
By continuity of the $K_0$ functor, we have
$$
(K_0(C^*_u(\Gamma)), K_0(C^*_u(\Gamma))^+) \cong \varinjlim ( \ell^\infty(\N,\Z), \ell^\infty(\N,\Z)^+, \alpha_n ).
$$
To describe the inductive limit more explicitly, define
$$
H_\Gamma^{(n)} := \left\lbrace  (m_1,m_2,...)\in \ell^\infty(\N,\Z): \sum_{i=jk_n+1}^{(j+1)k_n} m_i = 0 \text{ for all } j=0,1,2,...\right\rbrace
$$
and $H_\Gamma := \bigcup_{n=0}^\infty H_\Gamma^{(n)}$ (note that $\{H_\Gamma^{(n)}\}_{n=0}^\infty$ is an increasing sequence of subgroups of $\ell^\infty(\N,\Z)$ since $k_n$ divides $k_{n+1}$).

\begin{prop} Let $\alpha_n$ and $H_\Gamma$ be as above. Then
	$$
	\varinjlim(\ell^\infty(\N,\Z), \ell^\infty(\N,\Z)^+, \alpha_n) \cong (\ell^\infty(\N,\Z)/H_\Gamma, \ell^\infty(\N,\Z)^+/H_\Gamma),
	$$
	where $\ell^\infty(\N,\Z)^+/H_\Gamma$ is the collction of all elements in $\ell^\infty(\N,\Z)/H_\Gamma$ which can be represented by positive sequences.
\end{prop}
\begin{proof}
	First let us recall the standard construction of the inductive limit of ordered abelian groups: given an inductive system $\{ G_n, \alpha_n \}_{n=0}^\infty$ of ordered abelian groups, let $\nu_k:G_k\to \prod_{n=0}^\infty G_n$ be the map
	$$
	\nu_n(g) = (0,0,...,0, g, \alpha_{k}(g), \alpha_{k+1}\circ\alpha_{k}(g),...),
	$$
	where $g$ is in the $k$-th position. Define
	$$
	\beta_k := \pi\circ \nu_k:G_k\to \prod_{n=0}^\infty G_n/ \bigoplus_{n=0}^\infty G_n,
	$$
	where $\pi$ is the quotient map. Then $\{ \beta_k(G_k) \}_k$ is an incresing sequence of subgroups of $\prod_{n=0}^\infty G_n/ \bigoplus_{n=0}^\infty G_n$, and the unions
	$$
	G:= \bigcup_{k=0}^\infty \beta_k(G_k),\;\;\;\;\;\; G^+ := \bigcup_{k=0}^\infty \beta_k(G_k^+).
	$$
	form the inductive limit. Observe that when each $\alpha_n$ is surjective and $\alpha_n(G_n^+) = G_{n+1}^+$ (as in our case), we have
	$$
	\beta_0(G_0) = \beta_1(G_1) = \cdots,\;\;\;\;\;\; \beta_0(G_0^+) = \beta_1(G_1^+) = \cdots.
	$$
	Therefore the inductive limit $(G,G^+)$ is isomorphic to $\left(  G_k/\ker(\beta_k), G_k^+/\ker(\beta_k) \right)$ for any $k\in \N$. 
	
	Apply the discussion above to our case and choose $k =0$. It remains to show that $\ker(\beta_0)$ is equal to $H_\Gamma$. This is a straighforward computation (though probably a little confusing in the first glance, as we are dealing with sequences of sequences). Let 
	$$
	\mathfrak{m} = (m_1,m_2,m_3,...)
	$$
	be an element in $\ell^\infty(\N,\Z)$. Then
	$$
	\beta_0(\mathfrak{m}) = \left[ (m_1,m_2,...), \left( \sum_{i=1}^{k_1}m_i, \sum_{i=k_1+1}^{2k_1}m_i,...  \right), \left( \sum_{i=1}^{k_2}m_i, \sum_{i=k_2+1}^{2k_2}m_i,...  \right),...   \right],
	$$
	where the square bracket denotes the equivalence class in $\frac{\prod_{n=0}^\infty\ell^\infty(\N,\Z)}{\bigoplus_{n=0}^\infty\ell^\infty(\N,\Z)}$. Now $\beta_0(\mathfrak{m})$ vanishes if and only if there exists some $n\in \N$ such that
	$$
	\left( \sum_{i=1}^{k_n}m_i, \sum_{i=k_n+1}^{2k_n}m_i,...  \right) = (0,0,...)
	$$
	in $\ell^\infty(\N,\Z)$. From this description it is readily seen that $\ker(\beta_0)$ is equal to $H_\Gamma$.
\end{proof}

Now we can explicitly write down the ordered $K_0$-group of $C^*_u(\Gamma)$ for any countable locally finite group $\Gamma$.

\begin{thm}\label{thm4.8}
	Let $\Gamma$ be a countable locally finite group, and let $\{e\} =: \Gamma_0\subseteq \Gamma_1\subseteq \cdots \Gamma$ be an incresing sequence of finite subgroups such that $\Gamma = \bigcup_{n=0}^\infty \Gamma_n$. Define $k_n := |\Gamma_n|$, 
	$$
	H_\Gamma^{(n)} := \left\lbrace  (m_1,m_2,...)\in \ell^\infty(\N,\Z): \sum_{i=jk_n+1}^{(j+1)k_n} m_i = 0 \text{ for all } j=0,1,2,...\right\rbrace,
	$$
	and
	$$
	H_\Gamma := \bigcup_{n=0}^\infty H_\Gamma^{(n)}.
	$$
	Then 
	$$
	( K_0(C^*_u(\Gamma)), K_0(C^*_u(\Gamma))^+, [1]_0  ) \cong \left(  \ell^\infty(\N,\Z)/H_\Gamma, \ell^\infty(\N,\Z)^+/H_\Gamma, [\mathbf{1}]  \right),
	$$
	where $\mathbf{1}$ is the constant sequence with value 1.
\end{thm}
\begin{proof}
We only need to keep track of the order unit $[1]_0$. the $K_0$-class of the unit of $\prod_{i=1}^\infty \mathbb{C}$ in $\ell^\infty(\N,\Z) =: G_0$ is given by the constant sequence 1. Now the result follows since the structure map $G_0\to G_0/H_\Gamma$ for the inductive limit is nothing but the quotient map. 
\end{proof}

\begin{example} (cf. Example \ref{example:Z_2})
	Let us consider again the case that $\Gamma = \bigoplus_{\N}\Z/2\Z$ with the increasing sequence of subgroups $\Gamma_n = \bigoplus_{i=1}^n \Z/2\Z$. Then 
	$$
	H_\Gamma^{(1)} = \left\lbrace  (m_1,m_2,...)\in \ell^\infty(\N,\Z): m_1+m_2 = m_3+m_4 = \cdots = 0\right\rbrace,
	$$
	$$
	H_\Gamma^{(2)} = \left\lbrace  (m_1,m_2,...)\in \ell^\infty(\N,\Z): m_1+m_2 + m_3+m_4 = m_5+m_6+m_7+m_8 = \cdots = 0
	\right\rbrace,
	$$
	and so on. 
	
	In partiuclar, in the quotient group $\ell^\infty(\N,\Z)/H_\Gamma$ the constant sequence $\mathbf{1}$ can also be represented by the sequence $(2,0,2,0,2,0,2,0,...)$ or the sequence $(4,0,0,0,4,0,0,0,...)$.
\end{example}

\begin{cor}
Let $\Gamma$ and $\Lambda$ be two countable locally finite groups. Then $\Gamma$ and $\Lambda$ are coarsely equivalent if and only if $K_0(C_u^*(\Gamma))\cong K_0(C_u^*(\Lambda))$.
\end{cor}
\begin{proof}
If $\Gamma$ and $\Lambda$ are coarsely equivalent, then by \cite[Corollary 3.6]{STY02} $C_u^*(\Gamma)$ and $C_u^*(\Lambda)$ are Morita equivalent (see also \cite[Theorem~4]{BNW07}). In particular, they have isomorphic $K_0$ groups.

Conversely, if $\Gamma$ and $\Lambda$ are not coarsely equivalent, then by Remark~\ref{BZ} it must be the case that one is finite and the other is infinite. Let us assume that $\Lambda$ is finite. Then we have $K_0(C^*_u(\Gamma)) \cong \ell^\infty(\N,\Z)/H_\Gamma$ by Theorem~\ref{thm4.8} and $K_0(C^*_u(\Lambda)) \cong \Z$. It is not hard to see that these two groups are not isomorphic. For instance, one can show that $\ell^\infty(\N,\Z)/H_\Gamma$ is not singly generated.
\end{proof}
Finally, we are ready to give a proof of our main theorem, which provides the first classification result for a class of non-separable unital $C^*$-algebras.

\begin{thm}[Theorem \ref{thm:main_intro}]\label{thm:main}
	Let $\Gamma$ and $\Lambda$ be countable locally finite groups with proper left-invariant metrics $d_\Gamma$ and $d_\Lambda$, respectively. Then the following are equivalent:
	\begin{enumerate}
		\item $(\Gamma, d_\Gamma)$ and $(\Lambda, d_\Lambda)$ are bijectively coarsely equivalent;
		\item there is a $*$-isomorphism $\varphi:C^*_u(\Gamma)\to C^*_u(\Lambda)$ such that $\varphi(\ell^\infty(\Gamma)) = \ell^\infty(\Lambda)$; 
		\item $C^*_u(\Gamma)$ and $C^*_u(\Lambda)$ are $*$-isomorphic;
		\item $(K_0(C^*_u(\Gamma)), K_0(C^*_u(\Gamma))^+, [1]_0 ) \cong (K_0(C^*_u(\Lambda)), K_0(C^*_u(\Lambda))^+, [1]_0 )$;
		\item $( K_0(C^*_u(\Gamma)), [1]_0  ) \cong ( K_0(C^*_u(\Lambda)), [1]_0   )$. 
	\end{enumerate}
\end{thm}
\begin{proof}
$(1)\Longrightarrow (2)$: See Proposition \ref{prop:bce}.

$(2)\Longrightarrow (3)\Longrightarrow (4)\Longrightarrow (5)$: Obvious.

$(5)\Longrightarrow (1)$: To shorten the notations we write $A := C^*_u(\Gamma)$ and $B:= C^*_u(\Lambda)$. Let
$$
\varphi:(K_0(A),[1_A]_0)\to ( K_0(B),[1_B]_0   )
$$
be an isomorphism. 

Assume for the contrary that $(\Gamma,d_\Gamma)$ and $(\Lambda, d_\Lambda)$ are not bijectively coarsely equivalent. Then by Theorem \ref{thm:supernatural} the associated supernatural numbers $s(\Gamma)$ and $s(\Lambda)$ are not equal. Without losing generality we may assume there exist a prime number $p$ and a positive integer $r$ such that $p^r$ divides $s(\Gamma)$ but not $s(\Lambda)$. Let $\{\Gamma_n\}_{n=0}^\infty$ and $\{\Lambda_n\}_{n=0}^\infty$ be increasing sequences of finite subgroups of $\Gamma$ and $\Lambda$, respectively, such that $\Gamma = \bigcup_{n=0}^\infty \Gamma_n$ and $\Lambda = \bigcup_{n=0}^\infty \Lambda_n$. Define $k_n := |\Gamma_n|$,   $H_\Gamma^{(n)}$, and $H_\Gamma$ in the same way as before. Similarly, we have $k'_n := |\Lambda_n|$, $H_\Lambda^{(n)}$, and $H_\Lambda$ coming from the group $\Lambda$.

Since $p^r$ divides $s(\Gamma)$, there exists an element $[q]_0$ in $K_0(A)$ such that $p^r([q]_0) = [1_A]_0$. Indeed, if $p^r$ divides $s(\Gamma)$ then by definition $p^r$ divides $|\Gamma_n|$ ($= k_n$) for some $n$ (in the usual sense). Now one can take $q$ to be any projection in $\prod_{i=1}^\infty M_{k_n}(\C)$ which has pointwise rank $k_n/p^r$.

Applying the isomorphism $\varphi$, we obtain an element $[q']_0 := \varphi( [q]_0 )$ in $K_0(B)$ such that $p^r([q']_0) = [1_B]_0$. Write
$$
[q']_0 = [(m_1,m_2,...)]\in \ell^\infty(\N,\Z)/H_\Lambda.
$$
Then the equality $p^r([q']_0) = [1_B]_0$ implies that
$$
(p^rm_1-1,p^rm_2-1,...) \in H_\Lambda^{(n)}
$$
for some positive integer $n$. By the definition of $H_\Lambda^{(n)}$, we have (among other things) 
$$
p^r(m_1+m_2+\cdots + m_{k_n'}) - k_n' = 0,
$$
which is impossible because $p^r$ does not divide $k_n'$ (otherwise $p^r$ would divide $s(\Lambda))$. This completes the proof.
\end{proof}

\begin{rem}\label{(1)=(2)}
It is shown in \cite[Corollary~2.21]{L17} that two countable discrete groups are bijectively coarsely equivalent if and only if their canonical actions on Stone-\v{C}ech compactifications are continuously orbit equivalent. Since the proof of \cite[Proposition~4.13]{MR2460017} goes through without any change for \'{e}tale Hausdorff locally compact topologically principal $\sigma$-compact groupoids, it follows from \cite[Theorem~1.2]{L2_17} that (1) and (2) are equivalent for all countable discrete groups. 
\end{rem}

\begin{rem}\label{rigidity}
We would like to thank Rufus Willett for letting us know the following fact, which can be deduced from (the proof) of \cite[Theorem~1.1]{MR1700742}: if $X$ and $Y$ are uniformly discrete, bounded geometry non-amenable metric spaces, then any coarse equivalence between $X$ and $Y$ is close to a bijective coarse equivalence. Using this result, one deduces from the main theorem of \cite{SW13} that $(3)\Longrightarrow (1)$ holds for the class of non-amenable exact countable groups (cf. \cite[Corollary~6.2]{SW13}). Our main theorem provides a small class of amenable groups for which $(3)\Longrightarrow (1)$ holds.
\end{rem}

\begin{rem}
It is known from \cite[Corollary~8]{BZ11} that all countably infinite locally finite groups are coarsely equivalent. Thus, there are many coarsely equivalent locally finite groups with non-isomorphic uniform Roe algebras. For instance, $C_u^*(\bigoplus_{\N} \Z_2)\ncong C_u^*(\bigoplus_\N \Z_3)$ (see Theorem~\ref{thm:supernatural}). In fact, two countable direct sums of finite prime cyclic groups are bijectively coarsely equivalent if and only if they are isomorphic as groups (see \cite[Corollary~5.6]{MR2340955}). We refer the reader to \cite[Section 6]{LR17} for relevant discussions.
\end{rem}

\begin{rem}
It is shown by Rufus Willett and the first named author in \cite[Corollary~1.7]{LW} that $K_0(C_u^*(\Gamma))=0$ for every non-amenable countable group $\Gamma$ with asymptotic dimension one. Hence, we cannot expect Theorem~\ref{thm:main} to hold for general countable groups with asymptotic dimension one. For instance, let us consider the free group $F_n$ on $n$ generators and the wreath product group $\Z_2\wr F_n$, which is also finitely generated. Since they are non-amenable countable groups with asymptotic dimension one, their uniform Roe algebras have trivial $K_0$ groups. However, they are not coarsely equivalent (or equivalently, they are not quasi-isometric) as being finitely presented is invariant under quasi-isometries \cite[Proposition~V.4]{MR1786869} and $\Z_2\wr F_n$ is not finitely presented (see \cite[Theorem~1]{MR0120269}).
\end{rem}

\begin{question}
Does Theorem~\ref{thm:main} hold for all bounded geometry metric spaces with asymptotic dimension zero?
\end{question}

It follows from \cite[Section~11]{BZ11} that Theorem~\ref{thm:main} holds for all bounded geometry isometrically homogeneous metric spaces with asymptotic dimension zero. Recall that a metric space $X$ is called \emph{isometrically homogeneous} if for any two points $x,y\in X$ there is a bijective isometry $f:X\rightarrow X$ such that $f(x)=y$.

As a side note, it follows from \cite[Corollary~6.4]{MR2340955} and Proposition~\ref{prop:subset-of-N} that there are uncountably many coarsely inequivalent asymptotically 0-dimensional subspaces of $\Z_{\geq 0}$. On the other hand, there are only two coarsely equivalent classes of countable locally finite groups (cf. Remark~\ref{BZ}).


\section{Structure of uniform Roe algebras of locally finite groups}
In this section, we show that the reduced crossed product $C^*$-algebra $\ell^\infty(\Gamma)\rtimes_r\Gamma$ is locally finite-dimensional for any (not necessarily countable) locally finite discrete group $\Gamma$. Note that the countable case follows from \cite[Theorem 8.5]{WZ10}, since a countable group is locally finite if and only if it has asymptotic dimension zero (see Theorem \ref{thm:char-loc-fin}).

\begin{defn}
A $C^*$-algebra $A$ is called \emph{locally finite-dimensional} (or \emph{local AF}) if for every $a_1,\ldots, a_n\in A$ and $\epsilon>0$, there exist a finite-dimensional $C^*$-subalgebra $B$ of $A$ and elements $b_1,\ldots,b_n\in B$ such that $||a_i-b_i||<\epsilon$ for $i=1,\ldots,n$.
\end{defn}

\begin{prop}\label{local AF}
Let $\Gamma$ be a locally finite discrete group. Then the reduced crossed product $C^*$-algebra $\ell^\infty(\Gamma)\rtimes_r\Gamma$ is locally finite-dimensional, where $\Gamma$ acts on $\ell^\infty(\Gamma)$ by the left-translation.
\end{prop}
\begin{proof}
Note that $\ell^\infty(\Gamma)\rtimes_r\Gamma$ is the inductive limit of $\ell^\infty(\Gamma)\rtimes_r\Gamma_i$, where each $\Gamma_i$ is a finitely generated subgroup of $\Gamma$. Since $\Gamma$ is locally finite, it suffices to show that $\ell^\infty(\Gamma)\rtimes_r\Lambda$ is a (not necessarily sequential) inductive limit of finite-dimensional $C^*$-algebras for each finite subgroup $\Lambda$ of $\Gamma$. 

Toward this end, let $\Lambda$ be any finite subgroup of $\Gamma$ and $\Gamma=\bigsqcup_{s\in J}\Lambda s$ be the partition of $\Gamma$ into right cosets of $\Lambda$ in $\Gamma$. We may identify $\ell^\infty(\Gamma)$ with $\prod_{s\in J}\ell^\infty(\Lambda s)$ as $\Lambda$-$C^*$-algebras via the product of the restriction maps. In the following, we show that $\prod_{s\in J}\ell^\infty(\Lambda s)$ is an inductive limit of finite-dimensional $\Lambda$-$C^*$-algebras.

Now the proof is essentially an equivariant version of the argument used in \cite[Lemma 8.4]{WZ10}. Let $\text{FP}(J)$ be the directed set of all finite partitions of $J$ ordered by refinement. More precisely, we write $\mathcal{P}\leq \mathcal{Q}$ if and only if $\mathcal{Q}$ refines $\mathcal{P}$ (i.e., every element of $\mathcal{Q}$ is a subset of some element of $\mathcal{P}$). For  each $\mathcal{P}$ in $\text{FP}(J)$, define 
$$
 A_\mathcal{P} := 
 \left\lbrace
\begin{tabular}{ c|lc  } 
\multirow{2}{9.8em}{$(f_s)_{s\in J}\in \prod_{s\in J}\ell^\infty (\Lambda s)$} & $f_s(gs) = f_t(g t) \text{ for all } g\in \Lambda    \text{ whenever } s,t$ \\ & $\text{belong to the same member of } \mathcal{P}$ \\ 
\end{tabular}
\right\rbrace.
$$
In other words, $A_\mathcal{P}$ consists of sequences that are ``constant'' on each member of the partition $\mathcal{P}$. Since $\mathcal{P}$ is a finite partition, each $A_\mathcal{P}$ is a finite-dimensional subalgebra. Moreover, if $\mathcal{Q}$ is a finite partition of $J$ which refines $\mathcal{P}$, then $A_\mathcal{Q}$ contains $A_\mathcal{P}$. Finally, given any element $f$ in $\prod_{s\in J} \ell^\infty(\Lambda s)$, one can always find a finite partition $\mathcal{P}$ of $J$ so that $A_\mathcal{P}$ almost contains $f$. We conclude that $\prod_{s\in J} \ell^\infty(\Lambda s) = \lim_{\mathcal{P}\in \mathrm{FP}(J)} A_\mathcal{P}$ (as $\Lambda$-$C^*$-algebras), where the connecting maps are nothing but inclusions.

Since all connecting maps are injective (or we can use the fact that $\Lambda$ is finite), we conclude that $\ell^\infty(\Gamma)\rtimes_r \Lambda=\lim_{\mathcal{P}\in \text{FP}(J)} A_{\mathcal{P}}\rtimes_r \Lambda$ (see e.g. \cite[Lemma 2.5]{MR2010742}). As each $A_{\mathcal{P}}\rtimes_r \Lambda$ is a finite-dimensional $C^*$-algebra for every $\mathcal{P}\in \text{FP}(J)$, the proof is complete.
\end{proof}
\begin{rem}
By an almost identical proof, we see that $\ell^\infty(X)\rtimes_r \Gamma$ is locally finite-dimensional if $\Gamma$ is a locally finite discrete group acting on a set $X$. The only difference is to consider the partition of  $X=\bigsqcup_{x\in J} (\Lambda. x)$ into its $\Lambda$-orbits instead of right $\Lambda$-cosets, where $\Lambda$ is any finite subgroup of $\Gamma$.
\end{rem}
As an easy consequence, we give a summary of equivalent $C^*$-properties of uniform Roe algebras coming from locally finite groups (we refer the reader to \cite{Bla06} for the relevant concepts involved in the next theorem):

\begin{cor} \label{cor:structure}
	Let $\Gamma$ be a discrete group. Then the following are equivalent:
	\begin{enumerate}
	
		\item $\ell^\infty(\Gamma)\rtimes_r\Gamma$ is locally finite-dimensional;
		\item $\ell^\infty(\Gamma)\rtimes_r\Gamma$ is quasidiagonal;
		\item $\ell^\infty(\Gamma)\rtimes_r\Gamma$ has stable rank one;
		\item $\ell^\infty(\Gamma)\rtimes_r\Gamma$ has cancellation;
		\item $\ell^\infty(\Gamma)\rtimes_r\Gamma$ is stably finite;
		\item $\ell^\infty(\Gamma)\rtimes_r\Gamma$ is finite;
		\item $\Gamma$ is locally finite.
	\end{enumerate}
\end{cor}
\begin{proof}
$(1)\Longrightarrow (2)$: Since finite-dimenional $C^*$-algebras are quasidiagonal, this follows from the local characterization of quasidiagonality (cf. \cite[Lemma 7.1.3]{BO08}) and an application of Arveson's Extension Theorem (see, for example, the proof of \cite[Proposition 7.1.9]{BO08}).

$(2)\Longrightarrow (5)$: See e.g. \cite[Proposition V.4.2.6]{Bla06}.

$(1)\Longrightarrow (3)$: It follows from the fact that every finite-dimensional $C^*$-algebra has stable rank one.

$(3)\Longrightarrow (4)$: See e.g. \cite[Proposition V.3.1.24]{Bla06}.

$(4)\Longrightarrow (5)$: For every $n\in \N$, \cite[Proposition V. 2.4.14]{Bla06} implies that every isometry in $M_n(\ell^\infty(\Gamma)\rtimes_r\Gamma)$ is unitary. Hence, $\ell^\infty(\Gamma)\rtimes_r\Gamma$ is stably finite.

$(5)\Longrightarrow (6)$: It is clear from the definitions.

$(6)\Longrightarrow (7)$: It follows from \cite[Proposition 2.5]{Sca17}.

$(7)\Longrightarrow (1)$: This is Proposition~\ref{local AF}.
\end{proof}

\begin{rem} \label{rem:structure_space}
Let $(X,d)$ be a metric space with bounded geometry. Wei showed in \cite{Wei11} that quasidiagonality, stable finiteness and finiteness of $C_u^*(X)$ are all equivalent to $X$ being a ``box space'' provided that $X$ is infinite. In particualr $X$ can have arbitrarily large asymptotic dimension. On the other hand, Rufus Willett and the first named author showed in \cite{LW}  that the conditions
\begin{enumerate}
	\item $C^*_u(X)$ is AF,
	\item $C^*_u(X)$ is locally finite-dimensional,
	\item $C^*_u(X)$ has stable rank one, and
	\item $C^*_u(X)$ has cancellation
\end{enumerate}
are all equivalent to $X$ having asymptotic dimension zero.
\end{rem}

\ \newline
{\bf Acknowledgments}. We would like to thank Rufus Willett for helpful and enlightening discussions on the subject.


\bibliographystyle{plain}
\addcontentsline{toc}{chapter}{Bibliography}
\bibliography{Biblio-Database}

\end{document}